\documentclass{article}
\def\E{\ensuremath{{\bf E}}}
\def\R{\ensuremath{{\bf R}}}
\def\Z{\ensuremath{{\bf Z}}}

\def\N{\ensuremath{{\bf N}}}

\def\lh{\ensuremath{\underline{h}}}
\def\uh{\ensuremath{\overline{h}}}
\def\lp{\left(}
\def\rp{\right)}
\def\dist{\mbox{\rm dist}}
\def\dim{\mbox{\rm dim}}
\def\hau{\ensuremath{\mathcal{H}}}
\newcommand{\floor}[1]{{\ensuremath{[#1]}}}
\newcommand{\fp}[1]{\ensuremath{\{#1\}}}
\newcommand{\bin}[2]{\ensuremath{{#1 \choose #2}}}
\newtheorem{Theorem}{Theorem}
\newtheorem{Proposition}{Proposition}
\newtheorem{Lemma}{Lemma}
\newtheorem{Corollary}{Corollary}
\newtheorem{Definition}{Definition}
\newtheorem{Notation}{Notation}
\newtheorem{Example}{Example}
\newtheorem{Remark}{Remark}
\newtheorem{Hypothesis}{Hypothesis}

\def\BEq{\begin{equation}}
\def\EEq{\end{equation}}
\def\BHyp{\begin{Hypothesis}\rm}
\def\EHyp{\end{Hypothesis}}
\def\BRem{\begin{Remark}\rm}
\def\ERem{\end{Remark}}
\def\BEx{\begin{Example}\rm}
\def\EEx{\end{Example}}
\def\BThe{\begin{Theorem}}
\def\EThe{\end{Theorem}}
\def\BNot{\begin{Notation}\rm}
\def\ENot{\end{Notation}}
\def\BDef{\begin{Definition}\rm}
\def\EDef{\end{Definition}}
\def\BPro{\begin{Proposition}}
\def\EPro{\end{Proposition}}
\def\BLem{\begin{Lemma}}
\def\ELem{\end{Lemma}}
\def\BCor{\begin{Corollary}}
\def\ECor{\end{Corollary}}
\def\BProof{\noindent{Proof.~}}
\def\EProof{\\}

\begin{document}
\author{ \small \sc Marianne Clausel\thanks{
       Address:    Laboratoire d'Analyse et   de
Math\'ematiques Appliqu\'ees,   UMR 8050 du CNRS, Universit\'e Paris Est,
             61 Avenue du G\'en\'eral de Gaulle,
             94010 Cr\'eteil Cedex,
            France.
 $\;\; $ Email: clausel@univ-paris12.fr   }  and Samuel Nicolay\thanks{
     Institut de Math\'ematique, Grande Traverse, 12, B\^atiment B37, B-4000 Li\`ege (Sart-Tilman), Belgium.
 $\;\; $ Email: S.Nicolay@ulg.ac.be  } }

\date{ { \small Universit\'e Paris Est and Universit\'e de Li\`ege}  }
\title{ {\large \bf  Wavelets techniques for pointwise anti-H\"olderian irregularity} }
\maketitle
{\noindent\bf Abstract:} In this paper, we introduce a notion of
weak pointwise H\"older regularity, starting from the definition
of the pointwise anti-H\"older irregularity. Using this concept, a
weak spectrum of singularities can be defined as for the usual
pointwise H\"older regularity. We build a class of wavelet series
satisfying the multifractal formalism and thus show the optimality
of the upper bound. We also show that the weak spectrum of
singularities is disconnected from the casual one (denoted here
strong spectrum of singularities) by exhibiting a multifractal
function made of Davenport series whose weak spectrum
differs from the strong one.

{\noindent\bf Keywords} Pointwise H\"older regularity, Wavelets, Spectrum of singularities, Multifractal formalism.

{\noindent\bf Mathematics Subject Classification} 26A16, 42C40.

\section{Introduction}
The concept of H\"olderian regularity has been introduced to study
nowhere differentiable functions (several examples are given in
\cite{jaf-nic:08,jaf-nic:09}). An archetype of such functions is maybe
the Weierstra{\ss} function
\[
W_H(x)=\sum_{n=0}^{+\infty} a^{-nH} \cos(2\pi a^n x)\quad (0<H<1)
\]
exhaustively studied by Hardy in \cite{har:16}. He proved that for
every $a> 1$, this function is nowhere differentiable. More
precisely, the function $W_H$ satisfies the two following
conditions on $[0,1]^2$,
\[
|W_H(y)-W_H(x)|\le C_1 |x-y|^H
\]
and
\[
\sup_{(u,v)\in[x,y]^2}|W_H(u)-W_H(v)|\ge C_2 |x-y|^H
\]
for two constants $C_1$, $C_2$ that do not depend on $x$ or $y$.
The first inequality gives the regularity of $W_H$, which is said
uniformly H\"{o}lder with exponent $H$ on $[0,1]$. The second one
reflects the irregularity of the function; in particular, $W_H$ is
nowhere differentiable. One says, following \cite{tri:92}, that $W_H$
is uniformly anti-H\"{o}lder with exponent $H$ on $(0,1)$.

An increasing interest has been paid to
 functions $f$ that are both uniform H\"{o}lder and uniformly anti-H\"{o}lder with exponent $H$,
 since these two properties ensure that the
box-counting dimension of the graph of $f$ is equal to $2-H$ (see
e.g.\ \cite{fal:90}). Canonical Weierstra{\ss} functions,
i.e.\ functions of the form
\[
f(x)=\sum\limits_{n}b^{-n\alpha}g(b^{n}x)
\]
where $g$ is $1$-periodic, $1<b<\infty$ and $0<\alpha<1$, have
been extensively studied from the irregularity point of view by
many authors (see
\cite{kmy:84,mw:86,pu:89,bh:99,bh:00,heurt:03,heurt:05,dt:06,jaf-nic:08}).
Other well-known examples of such functions are provided by sample
paths of Gaussian fields, generalizing the fractional Brownian
motion. Irregularity properties, such as law of the iterated
logarithm, are established using fine results concerning the
regularity of the local time of the studied fields (see
\cite{ber:72,gh:80,adl:81}). In particular, in this class of
examples are included the so-called index-$\alpha$ Gaussian fields
studied in  \cite{adl:81}, or more generally non locally
deterministic Gaussian fields (see e.g. \cite{ber:73,pitt:78}) and
strongly non locally deterministic Gaussian fields (see
\cite{xiao:97,xiao:05,xiao:07}).

In this paper, we focus on the pointwise anti-H\"{o}lderian
irregularity, which is the pointwise counterpart of the concept of
uniform anti-H\"{o}lderian irregularity. Our main goal is to
answer quite natural questions: Can we overstep the usual
framework of functions both uniform H\"{o}lder and uniformly
anti-H\"{o}lder? More precisely, can we give some explicit
examples of functions for which the pointwise anti-H\"{o}lderian
behavior is different from point to point? What are the main characteristics of such a behavior?

In the case of the usual pointwise
regularity, multifractal functions provide examples of functions
for which the H\"{o}lder exponent vary from point to point. So, we
naturally tend to be interested in defining some multifractal
functions for this notion of pointwise anti-H\"{o}lderian
irregularity. We also need suitable tools to describe the
multifractal behavior of such functions. It raises the problem of
the related multifractal formalism. Indeed, let us recall that, in
general settings, it is not possible to estimate the regularity
index (which will be defined hereafter) of a function at a given
point. The relevant information is then the "size" of the sets of
points where the regularity is the same. This "size" is
mathematically formalized as the Hausdorff dimension. The function
that associates the dimension of the set of points sharing the
same regularity index with this index is referred to as the
spectrum of singularities. The goal of any multifractal formalism
is to provide a method which allows to estimate this spectrum of
singularities from numerically computable quantities derived from
the signal. The same problem arises when dealing with pointwise
anti-H\"{o}lderian irregularity.

Section~\ref{sec2} is devoted to the definitions related to the
H\"older regularity. In Section~\ref{sec:expholdsup} we investigate
the structure of the irregularity exponent and define, by means of
wavelet series, functions with prescribed irregularity exponent. In
Section~\ref{sec3}, we recall already known results about the
multifractal formalism for the pointwise anti-H\"{o}lderian
irregularity. Section \ref{sec4} is devoted to the question of the
validity of this multifractal formalism: Using multifractal
measures, we define a class of wavelet series for which the
multifractal formalism holds. In the last section, we compare the
two concepts of multifractal functions: The usual one and this new
one related to anti-H\"{o}lderianity. We show that the two notions
are clearly disconnected. Indeed, we exhibit an example of
Davenport series which is multifractal for the usual pointwise regularity but
monofractal for the pointwise irregularity.

\section{Pointwise H\"olderian regularity}\label{sec2}

We start by giving the definitions of the pointwise H\"{o}lderian
regularity and anti-H\"{o}lderian irregularity. The concept of
anti-H\"{o}lderian functions with exponent $H$ has been introduced
by C.Tricot in \cite{tri:92}; he formalized a notion already
used for investigating Weierstra{\ss}-type functions or sample paths properties
of locally non deterministic
Gaussian fields. Anti-H\"{o}lderian functions with
exponent $H$ were only defined in the case $H\in (0,1)$. A
consistent definition is given here for $H$ larger than
$1$. Since the anti-H\"olderian condition is stronger than just
negating the H\"olderian condition, a weaker H\"olderian
regularity is obtained by negating the anti-H\"olderian condition.
Finally, discrete wavelet transform and multiresolution analysis
are particularly efficient tools to study the H\"olderian regularity
of a function (see e.g.\ \cite{jaf:04}). The main results binding the
regularity of a function and its wavelet coefficients are briefly reviewed
at the end of this section.

Let us point out that the anti-H\"{o}lderian irregularity
condition has also been considered in the measure setting (see e.g.\ \cite{bmp:92}); a review of this measure-based irregularity framework is presented in Section \ref{subseq:Some results about multifractal analysis of measures}.

\subsection{Weak and strong pointwise H\"olderian regularity}
We recall first the definition of the H\"olderian regularity;
this definition naturally leads to a notion of H\"olderian irregularity. One will talk about H\"olderian and anti-H\"olderian functions. Finally, a weaker definition of pointwise smoothness is obtained by negating the condition related to the anti-H\"olderian functions.

\BDef Let $f:\R^d\to\R^{d'}$ be a locally bounded function, let
$x_0\in\R^d$ and $\alpha\ge 0$; $f\in C^\alpha(x_0)$ if there
exist $C,R>0$ and a polynomial $P$ of degree less than $\alpha$ such that
\begin{equation}\label{eq:strong holder}
\|f(x)-P(x)\|_{L^\infty\lp B(x_0,r)\rp }\le C r^\alpha, \quad \forall r\le
 R.
\end{equation}
Such a function is said H\"olderian of exponent $\alpha$ at $x_0$. The lower H\"
older exponent of $f$ at $x_0$ is
\[
\lh_f(x_0)=\sup\{\alpha:f\in C^\alpha(x_0)\}.
\]
A function $f$ is uniformly H\"olderian of exponent $\alpha$ ($f\in
C^\alpha(\R^d)$) if there exists $C>0$ such that (\ref{eq:strong
holder}) is satisfied for any $x_0\in\R^d$ and $R=\infty$; $f$ is
uniformly H\"olderian if there exists $\varepsilon>0$ such that $f\in
C^\varepsilon(\R^d)$. \EDef
Recall that the lower H\"older exponent
is simply denoted H\"older exponent in the literature. Since we are interested in introducing another concept of
pointwise H\"olderian regularity, the accustomed notation $h$ is replaced here by $
\lh$

The irregularity of a function can be studied through the notion of anti-H\"olderianity. Recall that the finite differences of arbitrary order are defined as follows,
\[
\Delta^{\! 1}_h f(x)=f(x+h)-f(x), \quad \Delta^{\! n+1}_h f(x)=\Delta^{\! n}_h f
(x+h)- \Delta^{\! n}_h f(x).
\]
We use the following notation,
\[
B_h(x_0,r)=\{x:[x,x+(\floor{\alpha}+1)h]\subset B(x_0,r)\}.
\]
Since condition (\ref{eq:strong holder}) is equivalent to
\begin{equation}\label{eq:strong holder2}
\sup_{|h|\le r} \| \Delta^{\! \floor{\alpha}+1}_h f \|_{L^\infty\lp B_h(x_0,r)\rp}\le C r^\alpha, \quad \forall r\le R
\end{equation}
(see e.g.\ \cite{dev:84,clau:08,kra:83}), the next definition is
intuitive. \BDef Let $f:\R^d\to\R^{d'}$ be a locally bounded
function, let $x_0\in\R^d$ and $\alpha\ge 0$; $f\in I^\alpha(x_0)$
if there exist $C,R>0$ such that
\begin{equation}\label{eq:anti-holder}
\sup_{|h|\le r} \|\Delta^{\!\floor{\alpha}+1}_h f \|_{L^\infty\lp B_h(x_0,r)\rp
}\ge C r^\alpha, \quad \forall r\le R.
\end{equation}
Such a function is said anti-H\"olderian of exponent $\alpha$ at $x_0$. Let us notice that the Whitney theorem asserts that $I^{\alpha_1}(x_0)\subset I^{\alpha_2}(x_0)$ if $\alpha_1\le \alpha_2$ (more precisely, it is a direct consequence of Proposition~1 of \cite{clanicdeux:09}). The upper H\"older exponent (or irregularity exponent) of $f$ at $x_0$ is
\[
\uh_f(x_0)=\inf\{\alpha:f\in I^\alpha(x_0)\}.
\]
We will say that $f$ is strongly H\"olderian of exponent
$\alpha$ at $x_0$ ($f\in C_s^\alpha(x_0)$) if $f\in
C^\alpha(x_0)\cap I^\alpha(x_0)$. \EDef

It follows from the definitions that if a function $f$ is anti-H\"olderian with exponent $\alpha$, it cannot be H\"olderian with exponent $\beta$ if $\beta>\alpha$. We thus have the following relation between the lower and upper exponents of $f$: $\lh_f \le \uh_f$.

The statement (\ref{eq:anti-holder}) is stronger than just negating
the H\"olderian regularity since such a negation only yields the
existence, for any $C>0$, of a subsequence $(r_n)_n$ (depending on
$C$) for which
\[
\sup_{|h|\le r_n} \|\Delta^{\!\floor{\alpha}+1}_h
f \|_{L^\infty\lp B_h(x_0,r_n)\rp }\ge C r_n^\alpha.
\]
We are naturally led to the following
definition.
\BDef Let $f:\R^d\to\R^{d'}$ be a locally bounded
function, let $x_0\in\R^d$ and $\alpha\ge 0$; $f\in
C^\alpha_w(x_0)$ if $f\notin I^\alpha(x_0)$, i.e.\ for any $C>0$
there exists a decreasing sequence $(r_n)_n$ such that
\begin{equation}\label{eq:Cw}
\sup_{|h|\le r_n} \|\Delta^{\!\floor{\alpha}+1}_h f \|_{L^\infty\lp B_h(x_0,r_n)\rp }\le C r_n^\alpha, \quad \forall n\in\N.
\end{equation}
Such a function is said weakly H\"olderian of exponent $\alpha$ at
$x_0$. \EDef Roughly speaking, a function is weakly H\"olderian of
exponent $\alpha$ at $x_0$ if for any $C>0$, one can bound the
oscillation of $f$ over $B(x_0,r_n)$ by $C r_{n}^{\alpha}$ for a
remarkable decreasing subsequence $(r_n)_n$ of scales, whereas for an
 H\"{o}lderian function, the oscillation of $f$ over
$B(x_0,r)$ has to be bounded at each scale $r>0$ by $C
r^{\alpha}$, for some $C>0$.

\subsection{H\"olderian regularity and wavelet coefficients}
Here, we review the wavelet criterion for strong H\"olderian regularity
and irregularity.

Let us briefly recall some definitions and notations (for more
precisions, see e.g.\ \cite{dau:88,mey:90,mal:98}). Under some
general assumptions, there exists a function $\phi$ and $2^d-1$
functions $(\psi^{(i)})_{1\le i<2^d}$, called wavelets, such that
 $\{\phi(x-k)\}_{k\in\Z^d}\cup\{\psi^{(i)}(2^j x-k):1\le
i<2^d, k\in \Z^d, j\in\Z \}$ form an orthogonal basis of
$L^2(\R^d)$. Any function $f\in L^2(\R^d)$ can be decomposed as
follows,
\[
f(x)=\sum_{k\in \Z^d} C_k \phi(x-k) + \sum_{j=1}^{+\infty}
\sum_{k\in\Z^d} \sum_{1\le i<2^d} c^{(i)}_{j,k} \psi^{(i)}(2^j
x-k),
\]
where
\[
c^{(i)}_{j,k}=2^{dj}\int_{\R^d}f(x) \psi^{(i)}(2^jx-k)\, dx,
\]
and
\[
C_k=\int_{\R^d} f(x) \phi(x-k)\, dx.
\]
Let us remark that we do not choose the $L^2$ normalization for
the wavelets, but rather an $L^\infty$ normalization, which is
better fitted to the study of the H\"olderian regularity. Hereafter, the
wavelets are always supposed to belong to $C^r$ with
$r>\alpha$ and the functions
$\{\partial^{s}\phi\}_{|s|\le r}$,
$\{\partial^{s}\psi^{(i)}\}_{|s|\le r}$ are assumed to have
fast decay.

A dyadic cube of scale $j$ is a cube of the form
\[
\lambda=[\frac{k_1}{2^j},\frac{k_1+1}{2^j})\times \cdots
\times [\frac{k_d}{2^j},\frac{k_d+1}{2^j}),
\]
where $k=(k_1,\ldots,k_d)\in \Z^d$. In the sequel, we denote $|\lambda|$ the scale of a dyadic cube $|\lambda|$.
From now on, wavelets and wavelet coefficients will be indexed
with dyadic cubes $\lambda$. Since $i$ takes $2^d-1$ values, we can assume
that it takes values in $\{0,1\}^d-(0,\ldots,0)$; we will use the
following notations:
\begin{itemize}
\item $\lambda=\lambda(i,j,k)=\frac{k}{2^j}+\frac{i}{2^{j+1}}+[0,\frac{1}{2^{j+1}})^d$,
\item $c_\lambda=c^{(i)}_{j,k}$,
\item $\psi_\lambda=\psi^{(i)}(2^j x-k)$,
\item $e_\lambda=k/2^j$.
\end{itemize}
The pointwise H\"olderian regularity of a function is closely related
to the decay rate of its wavelet leaders. \BDef The wavelet
leaders are defined by
\[
d_\lambda=\sup_{\lambda'\subset\lambda} |c_\lambda|.
\]
Two dyadic cubes $\lambda$ and $\lambda'$ are adjacent if they are
at the same scale and if $\dist(\lambda,\lambda')=0$. We denote by
$3\lambda$ the set of $3^d$ dyadic cubes adjacent to $\lambda$ and
by $\lambda_j(x_0)$ the dyadic cube of side $2^{-j}$ containing
$x_0$. Then
\[
d_j(x_0)=\sup_{\lambda\subset 3\lambda_j(x_0)} d_{\lambda}.
\]
\EDef
The following theorem (Theorem 1 of \cite{jaf:04}) allows to ``nearly'' characterize the H\"olderian regularity by a decay condition
 on $d_j$ as $j$ goes to infinity.
\BThe\label{thm:wav-hol}
Let $\alpha>0$; if $f\in C^\alpha(x_0)$, then there exists $C>0$ such that
\begin{equation}\label{eq:s h:cond}
d_j(x_0)\le C2^{-\alpha j}, \quad \forall j\ge 0.
\end{equation}
Conversely, if (\ref{eq:s h:cond}) holds and if $f$ is uniformly H\"older, then
there exist $C,R>0$ and a polynomial $P$ of degree less than $\alpha$ such that
\[
\|f(x)-P(x)\|_{L^\infty\lp B(x_0,r)\rp} \le C r^\alpha \log\frac{1}{r},
\quad \forall r\le R.
\]
\EThe

To give necessary and sufficient conditions concerning the
irregularity, we suppose that the wavelets are compactly
supported and belong to $C^{[\alpha]+1}(\R^d)$; such wavelets are constructed in \cite{dau:88}. The result relies on the following lemma.
\BLem\label{lem:deuxml}
Let $f\in L^\infty_{\rm loc}(\R^d)$; the two following assertions are then equivalent:
\begin{enumerate}
\item there exists some $\beta>1$ such that, for any $C>0$, there exists a non decreasing sequence of integers $(j_n)$ such that
\begin{equation}\label{eq:cdeuxmicrohyp}
\sup_{j\le j_n} \{2^{j([\alpha]+1)} d_j(x_0)\} \le C \frac{2^{j_n([\alpha]+1-\alpha)}}{j_n^\beta},
\end{equation}
\item there exists some $\beta>1$ such that, for any $C>0$, there exists a strictly increasing sequence of integers $(j_n)$ such that, for any $\lambda$,
\begin{equation}\label{eq:cdeuxmicro}
|c_\lambda|\le C (\theta(2^{-|\lambda|})+\theta(|x_0-e_\lambda|)),
\end{equation}
where $\theta$ is a non decreasing function such that, if $j\in\{j_n,\ldots, j_{n+1}-1\}$ for some $n\in\N$,
\[
\theta(2^{-j})= \inf(\frac{2^{j_{n+1}([\alpha]+1-\alpha)} 2^{-j([\alpha]+1)}}{j_{n+1}^\beta},
\frac{2^{-j_n\alpha}}{j_n^\beta}).
\]
\end{enumerate}
\ELem
\BProof
Let us suppose that Property~(\ref{eq:cdeuxmicro}) holds. For any $j$, we have $\theta(2\cdot 2^{-j})\le 2^{[\alpha]+1}\theta(2^{-j})$. Moreover, if $\lambda'\subset 3\lambda_j (x_0)$, one has $|x_{0}-e_{\lambda'}|\le 4d 2^{-|\lambda'|}=2^{-(|\lambda'|-2-\log_2 d)}$. Therefore,
\begin{eqnarray*}
 2^{j([\alpha]+1)}d_j(x_0) &=& 2^{j([\alpha]+1)}\sup_{\lambda'\subset 3\lambda_j(x_0)}|c_{\lambda'}|\\
 &\le& C2^{j([\alpha]+1)}\sup_{\lambda'\subset 3\lambda_j(x_0)}
 \theta(2^{-|\lambda'|})(1+2^{(2+\log_2 d)([\alpha]+1)})\\
&\leq& C 2^{j([\alpha]+1)}2^{[\alpha]+1}(1+2^{(2+\log_2 d)([\alpha]+1)})\theta(2^{-j}),
\end{eqnarray*}
Since  $2^{([\alpha]+1)j}\theta(2^{-j})\le 2^{j_n([\alpha]+1-\alpha)}/j_n^\beta$ for any $j\le j_n$, inequality~(\ref{eq:cdeuxmicrohyp}) is satisfied.

Conversely, if inequality~(\ref{eq:cdeuxmicrohyp}) holds, we have $d_{j_n}(x_0)\le C 2^{-j_n\alpha}/j_n^\beta$ and therefore, since the sequence $(d_j(x_0))$ is non increasing, $d_{j}(x_{0})\le C 2^{-j_n\alpha}/j_n^\beta$ for any $n\in\N$ and any $j\ge j_n$. Moreover, we have
\[
d_j(x_0) \le 2^{-j([\alpha]+1)}\sup_{j'\le j_{n+1}} \{2^{j'([\alpha]+1)}d_{j'}(x_0)\}
 \le C \frac{2^{j_{n+1}([\alpha]+1-\alpha)}}{j_n^\beta} 2^{-j([\alpha]+1)}
\]
for any $n\in \N$, and any $j\le j_{n+1}$. These relations imply that the inequality $d_{j}(x_{0})\le C\theta(2^{-j})$ is valid for any $n\in\N$ and any $j\in\{j_n\cdots,j_{n+1}-1\}$.
Let us now fix $\lambda'$ and set $j=\sup\{m:\lambda'\subset 3\lambda_m(x_{0})\}$. By definition of $j$, one has $|c_{\lambda'}|\leq d_j(x_0)\le C\theta(2^{-j})$.

If $j=|\lambda'|$ or if $j=|\lambda'|+1$, using the fact that $\theta$ is non decreasing, one gets
\[
|c_{\lambda'}|\le C\theta(2\cdot 2^{-j})\le C 2^{[\alpha]+1}\theta(2^{-|\lambda'|})
 \le C 2^{[\alpha]+1}(\theta(2^{-|\lambda'|})+\theta(|x_0-e_{\lambda'}|)).
\]
If $j<|\lambda'|$, $\theta(2^{-j})$ is larger than $\theta(2^{-|\lambda'|})$. Nevertheless, one has $2^{-j-1}\le |x_{0}-e_{\lambda'}|$ and thus
\begin{eqnarray*}
|c_{\lambda'}| &\le& d_j(x_0)\le C\theta(2^{-j})
 \le C\theta(|x_0-e_{\lambda'}|) \\
 &\le& C 2^{[\alpha]+1}(\theta(2^{-|\lambda'|})+\theta(|x_0-e_{\lambda'}|)).
\end{eqnarray*}
In any case, property~(\ref{eq:cdeuxmicro}) is recovered.
\EProof
This Lemma is indeed a two-microlocal characterization of Property~(\ref{eq:cdeuxmicrohyp}). We need it to prove the following wavelet characterization of pointwise irregularity:
\BThe\label{thm:wav-hol2}
Let $\alpha>0$ and $f\in L^\infty_{\rm loc}(\R^d)$;
\begin{enumerate}
\item if there exists $C>0$ such that
\begin{equation}\label{eq:cs}
d_j(x_0)\ge C 2^{-j\alpha},
\end{equation}
for any $j\ge 0$, then $f\in I^\alpha(x_0)$,
\item conversely, suppose that $f$ is
uniformly H\"older; if $f\in I^\alpha(x_0)$ then, for any $\beta>1$,
there exists $C>0$ such that
\begin{equation}\label{eq:cn}
\sup_{j'\le j} 2^{j'([\alpha]+1)}d_{j'}(x_0)\ge C
\frac{2^{j([\alpha]+1-\alpha)}}{j^{\beta}},
\end{equation}
for any $j\ge 0$.
\end{enumerate}
\EThe
\BProof
The first statement is a direct consequence of Theorem~4 of \cite{clanicdeux:09}. Let us prove the second statement by contrapositive. Assume that Property~(\ref{eq:cn}) does not hold, which is equivalent to assume that Property~(\ref{eq:cdeuxmicrohyp}) is satisfied. We use Lemma~\ref{lem:deuxml} to prove that inequality~(\ref{eq:Cw}) is satisfied for $r_n=C 2^{-j_n}$. Let $x\in \R^d$ such that $[x,x+([\alpha]+1)h]\subset B(x_0,2^{-j_n})$. We have
\[
\Delta^{\! M}_h f(x)=\sum_{k}
C_k \Delta^{\! M}_h \phi_k(x)+\sum_\lambda
c_\lambda \Delta^{M}_h\psi_\lambda(x)=\sum_{j\geq 0} \Delta_{h}^{M}f_j,
\]
where $f_0(x)=\sum_k C_k\phi_k(x)$ and $f_j(x)=\sum_{i,k} c_{j,k}^{(i)}\psi^{(i)}(2^j x-k)$ if $j\ge 1$. Since $f$ is assumed to be uniformly H\"older, there exists some $\varepsilon>0$ such that
\[
\max(\sup_k |C_k|,\sup_j (2^{j\varepsilon} \sup_\lambda |c_{\lambda}|))< \infty.
\]
Let $\alpha'>\alpha$ and define $J_n=\floor{\displaystyle\frac{\alpha'j_n}{\varepsilon}}+1$. The proof of Theorem~1 of \cite{jaf:04} yields the following inequality
\begin{eqnarray*}\label{eq:majo1}
\left|\sum_{j\geq J_n}\Delta^{[\alpha]+1}_{h}f_j(x)\right|\leq C 2^{-J_n\varepsilon_0}\leq C 2^{-j_n\alpha'}
\end{eqnarray*}
for $n$ sufficiently large.

Since the wavelets are assumed to be compactly supported, there exists $\ell_0$ such that $\mathrm{supp}(\psi^{(i)})\subset (-2^{\ell_0},2^{\ell_0})$. Let us give an upper bound of $\sum_{j=0}^{j_n+\ell_0}\Delta_h^{\! [\alpha]+1} f_j(x)$. If $|x_0-e_{\lambda}|\ge 2^{-j+\ell_0+1}$ then $|x-e_\lambda|\ge 2^{-j+\ell_0}$ and $\psi_{\lambda}(x)=0$. Therefore
\[
\sum_{j=0}^{j_n+\ell_0}\Delta_h^{\! [\alpha]+1}f_j(x)= \sum_{j=0}^{j_n+\ell_0}\sum_{\lambda,\,|x_0-e_\lambda|\le  2^{-j+\ell_0+1}} c_\lambda \Delta_h^{\! [\alpha]+1}\psi_\lambda(x).
\]
The regularity of the wavelets implies that $\sum_{j=0}^{j_n+\ell_0}\sum_{\lambda,\,|x_0-\mu_{\lambda}|\le  2^{-j+\ell_0+1}}c_{\lambda}\psi_{\lambda}(x)$ belongs to $C^{[\alpha]+1}(\R^d)$. Hence
\begin{eqnarray*}
\lefteqn{|\sum_{j=0}^{j_n+\ell_0} \sum_{\lambda,\,|x_0-e_\lambda|\le  2^{-j+\ell_0+1}}c_{\lambda} \Delta_h^{\! [\alpha]+1}\psi_{\lambda}(x)|}&&\\
&\le& |h|^{[\alpha]+1}\sum_{j=0}^{j_n+\ell_0}\|\sum_{\lambda,\,|x_0-e_{\lambda}|\le  2^{-j+\ell_0+1}}c_\lambda \psi_\lambda(x)\|_{C^{[\alpha]+1}(\R^d)}
\end{eqnarray*}
We now use Lemma~\ref{lem:deuxml} and the wavelet characterization of the spaces $C^{[\alpha]+1}(\R^d)$ to deduce that
\begin{eqnarray*}
\lefteqn{\|\sum_{\lambda,\,|x_0-e_\lambda|\le  2^{-j+\ell_0+1}}c_\lambda \psi_\lambda (x)\|_{C^{[\alpha]+1}(\R^d)}}&&\\
 &\le& \sup_{\lambda,\,|x_0-e_\lambda|\le  2^{-j+\ell_0+1}}(2^{j[\alpha]+1}|c_{\lambda}|)\\
 &\le& C 2^{j[\alpha]+1}(1+2^{(\ell_0+1)([\alpha]+1)})\theta(2^{-j})).
\end{eqnarray*}
This leads to the following upper bound,
\begin{eqnarray*}
|\sum_{j=0}^{j_n+\ell_0}\Delta_{\! h}^{[\alpha]+1} f_j(x)|&\le& |h|^{[\alpha]+1}C(1+2^{(\ell_0+1)([\alpha]+1)})\sum_{j=0}^{j_n} 2^{j([\alpha]+1)}\theta(2^{-j})\\
&=&|h|^{[\alpha]+1}C(1+2^{(\ell_0+1)([\alpha]+1)})2^{j_n ([\alpha]+1-\alpha)}\\
&\le& C(1+2^{(\ell_0+1)([\alpha]+1)})2^{-j_n \alpha}.
\end{eqnarray*}

Let us now give an upper bound of $\sum_{j=j_n+\ell_0}^{J_n}\Delta_{\! h}^{[\alpha]+1} f_j(x)$. In that case, let us remark that if $|x_0-e_\lambda|\ge 2^{-j_n+1}$ then $|x-e_\lambda|\ge 2^{-j+\ell_0+1}$ and $\psi_\lambda(x)=0$. We have
\begin{eqnarray*}
\lefteqn{|\sum\limits_{j=j_n+\ell_0}^{J_n}\Delta_{h}^{[\alpha]+1}f_j(x)|} &&\\
&=& |\sum_{j=j_n+\ell_0+1}\sum_{\lambda,\,|x_0-e_\lambda|\le  2^{-j_n+1}}c_\lambda \Delta_{\! h}^{[\alpha]+1}\psi_\lambda (x)|\\
&\le& ([\alpha]+1)\sum_{j=j_n+\ell_0+1}^{J_n}\sum_{\lambda,\,|x_0-e_\lambda|\le  2^{-j_n+1}}|c_\lambda|\sup_j \sup_{x\in\R^d}\sum_{|\lambda|=2^{-j}}|\psi_\lambda (x)|\\
&\le& C([\alpha]+1)\left(\sup_j \sup_{x\in\R^d}\sum_{|\lambda|=2^{-j}}|\psi_\lambda (x)|\right)\sum_{j=j_n+\ell_0+1}^{J_n}\sum_{\lambda,\,|x_0-e_\lambda|\le  2^{-j_n+1}} \theta(2^{-j_n})\\
&\le& C([\alpha]+1)\left(\sup_j \sup_{x\in\R^d}\sum_{|\lambda|=2^{-j}}|\psi_{\lambda}(x)|\right)J_n \frac{2^{-j_n\alpha}}{j_n^{\beta}}\\
&\le& C([\alpha]+1)\left(\sup_j \sup_{x\in\R^d}\sum_{|\lambda|=2^{-j}}|\psi_{\lambda}(x)|\right)2^{-j_n\alpha}
\end{eqnarray*}
for $n$ sufficiently large. Gathering these relations, we obtain $f\in C^{\alpha}_w(x_0)$.
\EProof

Note that we do not have a wavelet characterization of the property
$\overline{h}_f(x_0)=\alpha$. Indeed, it is proved in \cite{clau:08}
that, even up to a logarithmic correction, neither
condition (\ref{eq:cs}) is necessary, nor condition (\ref{eq:cn}) is
sufficient. Nevertheless, one can characterize the stronger
property $\lh_f(x_0)=\uh_f(x_0)=\alpha$ using
wavelets. Indeed, Theorems \ref{thm:wav-hol} and
\ref{thm:wav-hol2} lead to the following corollary that we will
use in the sequel.
\BCor\label{cor:specialcase}
Let $\alpha>0$ and
suppose that $f$ is uniformly H\"older. We have $\lh_f(x_0)=\uh_f(x_0)=\alpha$ if and only if
\[
\lim_{j\to\infty} \frac{\log d_j(x_0)}{-j\log 2}=\alpha.
\]
\ECor
\BProof
The first point of Theorem~\ref{thm:wav-hol2} implies that if $\lim_j \log d_j(x_0)/-j\log 2=\alpha$, we have $\lh_f(x_0)=\uh_f(x_0)=\alpha$. Let us prove the converse result. Assume that for any $\varepsilon>0$, we have $f\in C^{\alpha-\varepsilon}(x_0)\cap I^{\alpha+\varepsilon}(x_0)$. For any $\beta>1$, the preceding Theorems imply the existence of a constant $C>0$ such that $d_j(x_0)\le C2^{-j(\alpha-\varepsilon)}$ and $\sup_{j'\le j} 2^{j'([\alpha]+1)}d_{j'}(x_0)\ge 2^{j([\alpha]+1-\alpha-\varepsilon)} j^{-\beta}/C$, for any $j\ge 0$. Let $a>\varepsilon+2\varepsilon/([\alpha]+1-\alpha+\varepsilon)$ and $b=(\beta+\varepsilon)/([\alpha]+1-\alpha+\varepsilon)$. The previous relations lead to
\[
\sup_{j'\le j(1-a)-b\log_2 j}2^{j'([\alpha]+1)}d_{j'}(x_0)<
 C \frac{2^{j([\alpha]+1-\alpha-\varepsilon)}-j\varepsilon}{j^{\beta}}< C^{-1} \frac{2^{j([\alpha]+1-\alpha-\varepsilon)}}{j^{\beta}}
\]
and
\[
\sup_{j(1-a)-b \log_2 j\le j'\le j}2^{j'([\alpha]+1)} d_{j'}(x_0)
 \ge C^{-1} \frac{2^{j([\alpha]+1-\alpha-\varepsilon)}}{j^\beta}\ge C^{-1} \frac{2^{j([\alpha]+1-\alpha-\varepsilon)}}{j^\beta},
\]
for any $j\ge 2\log_2(C)/\varepsilon$ . Since the sequence $(d_j(x_0))$ is non increasing, we have
\[
2^{j([\alpha]+1)}d_{j(1-a)-b\log_2 j}(x_0)\ge
 C^{-1} \frac{2^{j([\alpha]+1-\alpha-\varepsilon)}}{j^\beta},
\]
for any $j\ge 2\log_2(C)/\varepsilon$. By setting $\ell=j(1-a)-b\log_2 j$, the preceding relation can be rewritten
\[
2^{\frac{1}{1-a}(\ell+2b\log_2 \ell)([\alpha]+1)}d_\ell(x_0)
 \ge C^{-1}\frac{2^{j([\alpha]+1-\alpha-\varepsilon)}}{j^\beta}\ge C^{-1}\frac{2^{\frac{\ell}{1-a}([\alpha]+1-\alpha-\varepsilon)}}{(\frac{\ell}{1-a})^\beta}.
\]
for any $\ell\ge 2\log_2(C)(1-a)/\varepsilon$.
Using the relation $d_j(x_0)\le C2^{-j(\alpha-\varepsilon)}$, we then obtain
\[
\alpha-\varepsilon\le\liminf_{j\to\infty} \frac{\log d_j(x_0)}{-j\log 2}\le\limsup_{j\to\infty} \frac{\log d_j(x_0)}{-j\log 2}<\alpha+\varepsilon.
\]
Since this inequality holds for any $\varepsilon>0$, the required result follows.
\EProof

\section{Construction of functions with prescribed lower and upper H\"older
exponents}\label{sec:expholdsup} In this Section, we investigate
in detail the structure of the irregularity exponent of a
continuous function. In Section~\ref{sec:WF}, we first prove,
considering a Weierstrass-type function, that it is possible to
construct a continuous function with prescribed pointwise
H\"{o}lder exponent $H$ provided that $H$ satisfies ``good
properties''. In Section~\ref{sec:ws-same h}  we focus on describing all the
functions which are both the classical H\"{o}lder exponent and
irregularity exponent of a continuous function. Finally, we study
the case where classical pointwise H\"{o}lder exponent and
irregularity exponent may differ. In this special case, we
give a sufficient condition and a necessary condition
for a couple of functions $(\underline{H},\overline{H})$
to be respectively the pointwise H\"{o}lder exponent and the
irregularity
exponent of a continuous function.

\subsection{A generic Weierstra{\ss} function with prescribed H\"older
exponents}\label{sec:WF}
In the same spirit as in \cite{dlm:98},
we consider the Weierstra{\ss}-type function
\begin{equation}\label{EqDefWeierstrass}
W(t)=\sum\limits_{j=0}^{+\infty}\lambda^{-jH(t)}\sin(2\pi\lambda^jt).
\end{equation}
and study its pointwise regularity.

\BPro\label{PropWeierstrass}
Let $H$ be a $\beta$-H\"{o}lderian function from $[0,1]$ to
$[a,b]\subset (0,1]$, satisfying $\sup_{t\in
[0,1]}|H(t)|<\beta$. If $W$ is a function of the form (\ref{EqDefWeierstrass}),
where $\lambda$ is an integer larger than $1$, then
\[
W\in C_s^{H(t)}(t)= C^{H(t)}(t)\cap I^{H(t)}(t),\quad \forall t\in [0,1].
\]
\EPro The proof of this proposition relies on the two following
lemma, analogous to Lemma $14$ and Proposition $15$ of
\cite{jaf-nic:08}.
\BLem\label{LemInterm1}
Let $\lambda>1$ and
$(f_j)_{j\in\N}$ a sequence of bounded and Lipschitz functions on
$\R$ for which there exists $C>0$ such that
\[
\|f_j\|_{\infty}+\|f'_j\|_{\infty}\leq C.
\]
The function $f(t)=\sum\limits_{j=0}^{+\infty}\lambda^{-jH(t)}f_j(\lambda^j
t)$ belongs to $C^{H(t)}(t)$, for any $t\in [0,1]$.
\ELem
\BProof Let $t\in [0,1]$. For any
$j$, $|f_j(t)-f_j(s)|$ is bounded by $C|s-t|$ or $C$. Let
$j_0(t)=[-\log |s-t|/(H(t)\log \lambda)]$. We have
\begin{eqnarray*}
|f(t)-f(s)|&\leq& C\sum_{j=0}^{j_0-1}|\lambda^j t-\lambda^j
s|\lambda^{-jH(t)}+C\sum_{j\geq j_0}\lambda^{-jH(t)}\\
&&+C\sum_{j=0}^{+\infty}|\lambda^{-jH(t)}-\lambda^{-jH(s)}|\\
&\leq& C|t-s|C \lambda^{j_0(1-H(t))}+C\lambda^{-j_0H(t)}\\
&&+C |\log \lambda||t-s|^\beta\sum_{j=0}^{+\infty}j\lambda^{-j a}\\
&\leq &C|t-s|^{H(t)}
\end{eqnarray*}
using the mean value theorem, $\beta$-H\"{o}lderianity of $H$ and the fact that
$H([0,1])\subset (0,1]$. \EProof
\BLem\label{LemInterm2} Let
$(f_j)_{j\in\N}$ be a sequence of $1$-periodic $C$-Lipschitz
functions from $\R$ to $\R$ and
\[
 f(t)=\sum_{j=0}^{\infty}\lambda^{-jH(t)}f_j(\lambda^j t),
\]
where $\lambda$ is an integer larger than $1$.
Assume that there exists $\ell
\in\{-\lambda-1,\cdots,\lambda-1\}$ and an integer $J$ such that
\[
D=\inf\limits_{j\geq J}|f_j(\frac{\ell}{\lambda})-f_j(0)|>0.
\]
If $\displaystyle\frac{C\ell}{\lambda(\lambda-1)}\leq D$ then
\[
f\in I^{H(t)}(t),\quad \forall t\in [0,1].
\]
\ELem
\BProof Since
\[
f(t)-f(s)=\sum\limits_{j=0}^{\infty}\lambda^{-jH(t)}(f_j(\lambda^j
t)-f_j(\lambda^j s))
+\sum\limits_{j=0}^{\infty}(\lambda^{-jH(t)}-\lambda^{-jH(s)})f_j(s)
\]
Proposition $15$ of \cite{jaf-nic:08} yields that for some $C_0>0$,
\[
|\sum\limits_{j=0}^{\infty}\lambda^{-jH(t)}(f_j(\lambda^j
t)-f_j(\lambda^j s))|\geq C_0 |t-s|^{H(t)}.
\]
Moreover
\[
|\sum\limits_{j=0}^{\infty}(\lambda^{-jH(t)}-\lambda^{-jH(s)})f_j(s)|\leq
|t-s|^{\beta}
\]
with $\beta>H(t)$. It provides the required conclusion.
\EProof
By applying Lemma~\ref{LemInterm1} and Lemma~\ref{LemInterm2} to
$f_j=sin(2\pi\cdot)$, $C=1$, $\ell=[\lambda/4]$ and
$D=|\sin(\ell/\lambda)|$, Proposition~\ref{PropWeierstrass} is
then straightforward.

\subsection{Wavelet series-defined functions with similar lower and upper H\"{o}lder exponents at any point}\label{sec:ws-same h}
The aim of this section is to prove a result analogous to Theorem $1$ of \cite{dlm:98}. Here we extend the results stated in \cite{dlm:98} since we give a characterization of functions which are both the lower and the upper exponent of a continuous function and thus satisfy a stronger property than in Theorem $1$ of \cite{dlm:98}.

\BThe\label{ThCarExpSim} Let $f$ a
 continuous nowhere differentiable function defined on $[0,1]$ with similar lower and upper H\"{o}lder exponent at any
point. There exists a sequence of continuous functions
$(H_j)_{j\in\N}$ from $[0,1]$ to $[0,1]$ such that
\[
\lh(t)=\uh(t)=\lim_{j\to \infty} H_j(t)\quad \forall t.
\]
Conversely, if $H$ is a function from $[0,1]$ to $[0,1]$ such that
\[
H(t)=\lim_{j\to \infty}H_j(t),
\]
where the $(H_j)_{j\in\N}$
are continuous functions, then there exists a continuous function
$f$ defined on $[0,1]$ such that
\[
\lh(t)=\uh(t)=H(t),\quad \forall t.
\]
\EThe

The first part of Theorem \ref{ThCarExpSim} is
straightforward. If one sets
\[
 \omega_{x}(t)=\sup_{|h|\leq t} |f(x+h)-f(x)|,
\]
the sequence $(H_j)_{j\in\N}$ defined by
\[
H_j(t)=\displaystyle\frac{\log(\omega_{t}(2^{-j})+2^{-2j})}{-j\log 2}
\]
satisfies the required conditions, since
\[
\lh(t)=\uh(t)=
\lim_{j\to \infty} \frac{\log \omega_{t}(2^{-j})}{-j\log 2}.
\]
Let us prove the converse assertion by means of wavelet series. We
will need the following wavelet criterium for the pointwise
regularity (see \cite{dlm:98,jaf:04}): \BPro\label{pro:WC} Let
$\alpha>0$ and assume that there exists $C>0$ such that for any
$j\in\N$,
\begin{equation}\label{eq:hypWC}
\sup_{k}|c_{j,k}|\leq C 2^{-\frac{j}{\log j}},
\end{equation}
and
\[
d_j(t)\leq C2^{-j\alpha}.
\]
Then for any $\varepsilon>0$, $f\in C^{\alpha-\varepsilon}(t)$.
\EPro \BRem This Proposition is a reformulation of Proposition $4$
of \cite{dlm:98} in terms of wavelet coefficients. Let us point
out that in \cite{dlm:98}, $f$ is expanded in the Schauder basis.
Then for any $J$ (Lemma 1 of \cite{dlm:98}),
\[
\sum\limits_{j\leq J}\sum\limits_k c_{j,k}\Lambda_{j,k}
\]
is a continuous piecewise affine function coinciding with $f$ on
dyadic numbers. In the case of wavelet basis, this property does
not hold. In order to prove Proposition \ref{pro:WC}, we thus need
to assume (\ref{eq:hypWC}) and use different arguments.
\ERem

\BProof Since $\sup_{k}|c_{j,k}|\leq
C2^{-\frac{j}{\log j}}$, the wavelet series converge uniformly on
any compact. Let $(f_j)_{j\ge -1}$ be defined as follows,
\[
f_{-1}(x)=\sum\limits_{k}C_k\varphi(x-k),\quad f_j(x)=\sum\limits_{k}c_{j,k}\psi_{j,k}(x)\quad (\forall
j\geq 0)
\]
As in \cite{jaf:04}, if $|\beta|\leq [\alpha]$,
the series $\sum_{j}\partial^{\beta}f_j$ converges absolutely. Now, for any $j\geq -1$, let us define
\[
P_{j,t}(x)=\sum_{\beta\leq
[\alpha]}\frac{(x-t)^{\beta}}{\beta !}\partial
f_j(t).
\]
If $j_0$ is the number such that
\[
2^{-j_0-1}\leq |x-t|\leq 2^{-j_0},
\]
and $j_1$ satisfies $2^{-\frac{j_1}{\log j_1}}\leq 2^{-\alpha
j_0}$, then as proved in \cite{jaf:04},
\[
\sum_{j\leq j_0}|f_j(x)-P_{j,t}(x)|\leq Cj_0 |x-t|^{\alpha}
\]
and
\[
\sum_{j\geq j_0}|P_{j,t}(x)|\leq C |x-t|^{\alpha}.
\]
Moreover,
\[
\sum_{j= j_0}^{j_1}|f_j(x)|\leq C j_1 |x-t|^{\alpha}.
\]
Since for any $\varepsilon$ and $j_0$ sufficiently large,
\[
j_1\leq C_\varepsilon 2^{\varepsilon j_0},
\]
one has
\[
\sum_{j= j_0}^{j_1}|f_j(x)|\leq C_\varepsilon
|x-t|^{\alpha-\varepsilon}.
\]
Since
\[
\sum_{j\geq j_1}|f_j(x)|\leq C \sum_{j\geq j_1} 2^{-\frac{j}{\log j}}
\leq C \log^2(j_1)2^{-\frac{j_1}{\log j_1}}
\leq |x-t|^{\alpha-\varepsilon},
\]
the proposition follows.
\EProof
We will also use a slightly modified
version of Lemma $2$ of \cite{dlm:98}:
\BPro Let $H$ a function
from $[0,1]$ to $[0,1]$ such that $H(t)=\lim_{j}H_j(t)$, where $(H_j)_{j\in\N}$ is a sequence of continuous functions. There exists a sequence $(P_j)_{j\in\N}$ of polynomials such that
\begin{equation}\label{EqPolyApprox}
\left\{
\begin{array}{l}
H(t)=\lim\limits_{j\to \infty} P_j(t),\quad \forall t\in [0,1],\\
\|P'_j\|_{\infty}\leq j,\quad \forall j\in\N.
\end{array}
\right.
\end{equation}
\EPro
We can now define a wavelet series with the desired properties:
\BPro Let $H$ be a
function from $[0,1]$ to $[0,1]$ such that
$H(t)=\lim_{j}H_j(t)$, where the $(H_j)_{j\in\N}$
is a sequence of continuous functions. Let $(P_j)_{j\in\N}$ be a sequence of
polynomials satisfying the relations (\ref{EqPolyApprox}). For any
$(j,k)\in\N\times\{0,\cdots,2^j-1\}$, set
\[
H_{j,k}=\max (\frac{1}{\log j},P_j(\frac{k}{2^j})).
\]
The function $f$ defined as
\begin{equation}\label{EqWS}
f(x)=\sum_{j\in\N}\sum_{k=0}^{2^j-1}2^{-jH_{j,k}}\psi_{j,k}(x).
\end{equation}
satisfies the following relations,
\[
\underline{h}(t)=\overline{h}(t)=H(t),\quad \forall t\in [0,1].
\]
\EPro
Let us look at a particular case.
\BRem If $H$ is a continuous function, the wavelet series
\[
\sum_{j\in\N}\sum_{k=0}^{2^j-1}2^{-jH_{j,k}}\psi_{j,k}(x),
\]
with $H_{j,k}=\max(1/\log j,H(2^{-j}k))$ has $H$ both as lower and
upper H\"{o}lder exponents.

\BProof
If $\lambda'=\lambda'(j',k')\subset 3\lambda_j(t)$, then
\[
|\frac{k'}{2^{j'}}-t|\leq C_1 2^{-j}.
\]
Since the sequence $(P_j)_{j\in\N}$ satisfies equalities
(\ref{EqPolyApprox}), for any $\varepsilon>0$, there exists an
integer $j_0$ such that for any $j\geq j_0$ and $\lambda' \subset
3\lambda_j(t)$,
\[
|H_{j',k'}-H(t)|\leq
|P_j'(t)-P_j'(\frac{k'}{2^{j'}})|+|H(t)-P_j'(t)|\leq
j2^{-j}+|H(t)-P_j'(t)|\leq \varepsilon.
\]
Then for any $j\geq j_0$,
\[
\max(2^{-\frac{j'}{\log j'}},2^{-j'\varepsilon}2^{-j'H(t)})\leq
|c_\lambda'|\leq
\max(2^{-\frac{j'}{\log j'}},2^{j'\varepsilon}2^{-j'H(t)}).
\]
We deduce that
\[
\lim_j \frac{\log d_j(t)}{-j\log 2}=H(t),\quad \forall t\in [0,1]
\]
and therefore
\[
\uh(t)=\lh(t)=H(t),\quad \forall t\in [0,1].
\]
\EProof
\ERem

\subsection{Wavelet series-defined functions with different lower and upper H\"{o}lder exponents}
Our main goal is to prove the following Theorem.
\BThe\label{ThCarExpSim2} Let $f$ a continuous nowhere
differentiable function defined on $[0,1]$. There exists a sequence of continuous functions
$(H_j)_{j\in\N}$ such that
\[
\lh(t)=\liminf_{j\to \infty} H_j(t),\quad
\uh(t)=\limsup_{j\to \infty}H_j(t),\quad
\forall t\in [0,1].
\]

Conversely, let $(\underline{H},\overline{H})$ a couple of
functions from $[0,1]$ to $[a,b]\subset (0,1)$ such that
\[
\underline{H}(t)=\liminf_{j\to \infty}H_j(t),\quad
\overline{H}(t)=\limsup_{j\to \infty}H_j(t),
\]
where the $(H_j)_{j\in\N}$ is a sequence of continuous functions. There exists a
uniform H\"{o}lder function $f$ from $[0,1]$ to $\R$ such that
\[
\lh(t)=\underline{H}(t)\leq
\overline{H}(t)=\uh(t),\quad
\forall t\in [0,1].
\]
\EThe
\BRem
The second assertion of Theorem~\ref{ThCarExpSim2} is
much weaker than the corresponding one of Theorem~\ref{ThCarExpSim}. In order to ensure the existence of a function $f$ with prescribed lower and upper exponents at any point, we need stronger assumptions on $\underline{H}$ and $\overline{H}$: indeed we assume that these functions take values in $[a,b]\subset (0,1)$ (and not in $(0,1)$).
\ERem
\BRem Let $\underline{H}$ and $\overline{H}$
be two continuous functions from $[0,1]$ to $[a,b]\subset (0,1)$ satisfying
\[
\forall t\in [0,1],\,\underline{H}(t)\leq\overline{H}(t).
\]
Define the sequence $(H_j)_{j\in\N}$ as
\[
H_{2j}=\underline{H},\quad
H_{2j+1}=\overline{H},\quad
\forall j\in\N.
\]
Theorem~\ref{ThCarExpSim2} provides the existence of a
function $f$ with lower and upper exponent $\underline{H}$ and
$\overline{H}$ at any point.\ERem
The proof of the direct part of
Theorem \ref{ThCarExpSim2} is exactly similar to this of Theorem
\ref{ThCarExpSim} and is left to the reader. In order to prove the
converse assertion we first need the following lemma.
\BLem Let
$(\underline{H},\overline{H})$ a couple of functions defined from $[0,1]$ to $[a,b]\subset (0,1)$
such that
\[
\underline{H}(t)=\liminf_{j\to \infty}H_j(t),\quad
\overline{H}(t)=\limsup_{j\to \infty}H_j(t),
\]
where the $H_j$ are continuous functions from $[0,1]$ to $[0,1]$. For any $(a',b')\in (0,1)^2$ such that $a'<a<b<b'$, there exists a sequence $(P_j)_{j\in\N}$ of polynomials from
$[0,1]$ to $[a',b']$ and for any $t\in [0,1]$, there exists a strictly increasing
sequence of integers $(j_n(t))_{n\in\N}$ depending on $t$ such that the three
following properties hold simultaneously
\begin{itemize}
 \item $\forall t\in [0,1]$,
\begin{equation}\label{EqPolyApprox3}
\underline{H}(t)=\liminf_{j\to \infty}P_j(t),\quad
\overline{H}(t)=\limsup_{j\to \infty}P_j(t),
\end{equation}
 \item $\forall j\in\N$,
\begin{equation}\label{EqPolyApprox5}
\|P'_j\|_{\infty}\leq j,
\end{equation}
 \item $\forall t\in [0,1]$ $\forall n\in\N$, $\forall j\in\{j_n(t),\cdots,j_{n+1}(t)-1\}$,
\begin{equation}\label{EqPolyApprox4}
jP_j(t)\geq\sup(j_n(t)(\overline{H}(t)-\varepsilon),j+j_{n+1}(t)(\overline{H}(t)-\varepsilon-1)).
\end{equation}
\end{itemize}
\ELem
\BProof Lemma $2$ of \cite{dlm:98} implies
that there exists a sequence of polynomials $(Q_\ell)_{j\in\N}$
such that Conditions (\ref{EqPolyApprox3}) and
(\ref{EqPolyApprox5}) both hold. Moreover in the construction of
\cite{dlm:98}, one may assume
\[
\|Q_\ell-H_\ell\|_{\infty}\leq \frac{1}{\ell}.
\]
Then, for $\ell$ sufficiently large, $Q_\ell$ is onto
$[a',b']\subset (0,1)$. Set $\beta_1=[b'/a']+1$,
$\beta_2=[(1-a')/(1-b')]+1$, $\beta=\beta_1\beta_2$ and define the
sequence $(P_j)_{j\in\N}$ as follows
\[
\forall \ell\in\N,\,\forall \beta^\ell+1\leq j\leq \beta^{\ell+1},\quad
P_{j}=Q_{\ell}.
\]
We now prove that the sequence $(P_j)_{j\in\N}$ satisfies the
required properties. Let $\varepsilon>0$, $t\in [0,1]$,
$(\ell_n(t))_{n\in\N}$ a sequence such that
\[
Q_{\ell_n(t)}\geq \overline{H}(t)-\varepsilon
\]
and set $j_n(t)=\beta^{\ell_n(t)}\beta_2$. For any integer $j$ we
distinguish three cases.
\begin{itemize}
\item If there exists $n\in \N$ such that $j_n\leq \ell\leq
\beta_1 j_n$, then $P_j=Q_{\ell_n}\geq \overline{H}-\varepsilon$.
Hence
\[
jP_j\geq j(\overline{H}-\varepsilon)\geq j_n
(\overline{H}-\varepsilon),
\]
and
\[
j(1-P_j)\leq j(1-\overline{H}+\varepsilon)\leq j\beta_1
(1-\overline{H}+\varepsilon) j_{n+1}(1-\overline{H}+\varepsilon).
\]
\item If there exists $n\in \N$ such that $\beta_1 j_n\leq j\leq
j_{n+1}/\beta_2$, then $P_j=Q_{\ell}$ with $\ell_{n}+1\leq
\ell\leq \ell_{n+1}-1$. Then
\[
jP_j\geq ja'\geq \displaystyle\frac{j}{\beta_1}\overline{H}\geq
j_n\overline{H},
\]
and
\[
j(1-P_j)\leq j(1-a')\leq j\beta_2 (1-\overline{H})\leq
j_{n+1}(1-\overline{H}).
\]
\item If there exists $n\in \N$ such that $j_{n+1}/\beta_2\leq
j\leq j_{n+1}$, then $P_j=Q_{j_{n+1}}\geq
\overline{H}-\varepsilon$. Hence
\[
jP_j\geq j(\overline{H}-\varepsilon)\geq
j_n(\overline{H}-\varepsilon),
\]
and
\[
j(1-P_j)\leq j(1-\overline{H}+\varepsilon)\leq
j_{n+1}(1-\overline{H}+\varepsilon).
\]
\end{itemize}
Thus, in any case we obtain the required property.
\EProof
Now we
prove the following proposition
\BPro
Let
$(\underline{H},\overline{H})$ a couple of functions from $[0,1]$
to $[a,b]\subset (0,1)$ such that
\[
\underline{H}(t)=\liminf_{j\to \infty}H_j(t),\quad
\overline{H}(t)=\limsup_{j\to \infty}H_j(t),
\]
where the $(H_j)_{j\in\N}$ is a sequence of continuous functions. Let $(P_j)_{j\in\N}$ a
sequence of polynomials satisfying Properties
(\ref{EqPolyApprox3}),(\ref{EqPolyApprox5}) and (\ref{EqPolyApprox4}) and consider the wavelet series defined by
\[
f(x)=\sum_{j\in\N}\sum_{k=0}^{2^j-1}2^{-jP_j(\frac{k}{2^j})}.
\]
Then for any $t\in [0,1]$,
\[
\lh(t)=\underline{H}(t)<\uh(t)=\overline{H}(t).
\]
\EPro
\BProof The assumption $P_j([0,1])\subset
[a',b']\subset(0,1)$ implies that $f$ is uniform H\"{o}lder. Since
for any $j$, $P_j$ satisfies Properties
(\ref{EqPolyApprox3}), (\ref{EqPolyApprox5}) and (\ref{EqPolyApprox4}), for any $\lambda\in 3\lambda_j(x_0)$, we have
\[
2^{-jP_j(t)}2^{-j^2 2^{-j}}\leq |c_\lambda|\leq 2^{-jP_j(t)}2^{j^2
2^{-j}}.
\]
Thus, for any $\varepsilon>0$, there exists $j_0$
sufficiently large such that
\[
\forall j\geq j_0,\,2^{-jP_j(t)}2^{-\varepsilon j}\leq
d_{j}(x_0)\leq 2^{-jP_j(t)}2^{\varepsilon j}.
\]
By definition of $\underline{H}$,
\[
\liminf_{j\to \infty}\frac{\log d_j(t)}{-\log 2}
=\liminf_{j\to \infty} P_j(t)=\underline{H}(t).
\]
Hence,
\[
\lh(t)=\underline{H}(t).
\]
In the same way,
\[
\uh(t)\leq\overline{H}(t).
\]
We now use Properties
(\ref{EqPolyApprox3}), (\ref{EqPolyApprox5}) and (\ref{EqPolyApprox4}).
There exists a strictly increasing sequence of integers
$(j_{n})_{n\in\N}$ such that
(\ref{EqPolyApprox3}), (\ref{EqPolyApprox5}) and (\ref{EqPolyApprox4})
hold. Then, $\forall n\in\N$, $\forall j\in\{j_n,\cdots,j_{n+1}\}$,
\[
d_{j}(t)\leq
2^{-jP_j(t)}2^{\varepsilon j}\leq
\inf(2^{-j_{n}(\overline{H}(t)-\varepsilon)},2^{-j_{n+1}(\overline{H}(t)-\varepsilon-1)-j}).
\]
The wavelet criteria then provides
\[
\overline{h}(t)\geq\overline{H}(t).
\]
\EProof

\section{Weak multifractal formalism}\label{sec3}
The aim of the multifractal analysis is to study ``irregularly
irregular'' functions, i.e.\ functions whose H\"older exponent can
jump from point to point. From a practical point of view, the
numerical computation of the pointwise H\"older exponent of a
signal is completely instable, and is indeed quite meaningless,
especially for signals whose pointwise H\"older exponent can take
very different values. Leaving this utopian view, one rather
wishes to get global informations about the pointwise regularity:
What are the values taken by the H\"older exponent? What is the
``size'' of the set of points $E_h$ where the H\"older exponent
takes a given value $h$? First of all, one has to define this
notion of size. Since the sets under consideration can be dense or
negligible, by ``size'', we cannot mean ``Lebesgue measure''. The
``fractal dimensions'' are more fitted for this purpose. Once the
right definition of dimension has been chosen, one still has to
determine the spectrum of singularities of the function, i.e.
the dimension of the sets $E_h$. This is the purpose of the
multifractal formalism. Naturally, all these definitions can be
transposed for the upper H\"older exponent.

\subsection{A notion of dimension}
In multifractal analysis, the notion of dimension which is mainly
used is the Hausdorff dimension. We recall here its definition.

The Hausdorff dimension is defined through the Hausdorff measure (see
\cite{fal:90} for more details). The best covering of a set $E\subset \R^d$ with sets subordinated to a diameter $\varepsilon$ can be estimated as follows,
\[
\hau_\varepsilon^\delta(E)=\inf\{\sum_{i=1}^\infty |E_i|^\delta:E\subset \bigcup_{i=1}^\infty E_i, |E_i|\le \varepsilon\},
\]
where for any $i$, $|E_i|$ denotes the diameter of $E_i$.\\
Clearly, $\hau_\varepsilon^\delta$ is an outer measure. The Hausdorff measure is defined from $\hau_\varepsilon^\delta$ as $\varepsilon$ goes to $0$.
\BDef
The outer measure $\hau^\delta$ defined as
\[
\hau^\delta(E)=\sup_{\varepsilon>0} \hau_\varepsilon^\delta(E)
\]
is a metric outer measure. Its restriction to the $\sigma$-algebra
of the $\hau^\delta$-measurable sets defines the Hausdorff measure
of dimension $\delta$. \EDef Since the outer measure $\hau^\delta$
is metric, the algebra includes the Borelian sets. The Hausdorff
measure is decreasing as $\delta$ goes to infinity. Moreover,
$\hau^\delta(E)>0$ implies $\hau^{\delta'}(E)=\infty$ if
$\delta'<\delta$. The following definition is thus meaningful.
\BDef The Hausdorff dimension $\dim_\hau(E)$ of a set $E\subset
\R^d$ is defined as follows
\[
\dim_\hau(E)=\sup\{\delta:\hau^\delta(E)=\infty\}.
\]
\EDef
With this definition, $\dim_\hau(\emptyset)=-\infty$.

\subsection{From the strong multifractal formalism to the weak multifractal formalism}
We first review the wavelet leaders based multifractal formalism as defined by Jaffard \cite{jaf:04},
which is one of the two methods allowing to recover, in some particular cases, the entire spectrum of singularities
(the second one is the wavelet transform of the maxima of the modulus method, introduced by Arneodo and his collaborators \cite{arn-al:95}).
Other multifractal formalisms only give, at best, the increasing part of the spectrum (see \cite{jaf:97}). These considerations on the strong H\"olderian regularity can be transposed to the weak one.

The lower spectrum of singularities allows to characterize
globally the regularity of a function through the lower H\"older
exponents. \BDef Let $f$ be a locally bounded function; its lower
isoH\"older sets are the sets
\[
\underline{E}_H=\{x:\lh(x)=H\}.
\]
The lower spectrum of singularities of $f$ is the function
\[
\underline{d}: \R^+ \cup \{\infty\} \to \R^+ \cup \{-\infty\}
\quad H\mapsto\dim_\hau (\underline{E}_H),
\]
\EDef

It is not always possible to compute the lower spectrum of
singularities of a function. A multifractal formalism is a method
that is expected to yield the function $\underline{d}$ through the use of a
Legendre transform. These formalisms are variants of a seminal
derivation which was proposed by Parisi and Frisch
\cite{par-fri:85}. The wavelet leaders method (WLM) uses wavelet
coefficients instead of $L^p$ norms, which are meaningless for
negative values of $p$. The partition function is defined as
follows
\[
S(j,p)=2^{-j}\sum_{\lambda:|\lambda|=2^{-dj} } d_\lambda^p.
\]
By setting,
\begin{equation}\label{eq:omega}
\omega(p)=\liminf_{j\to\infty} \frac{\log S(j,p)}{\log 2^{-j}},
\end{equation}
the spectrum of singularities $\underline{d}(h)$ is expected to be equal to
\begin{equation}\label{eq:smf:leg}
\inf_p\{h p-\omega(p)+d\}.
\end{equation}

The heuristic argument leading to the previous method is the
following. The contribution of the dyadic cubes of side $2^{-j}$
containing a point whose lower H\"older exponent is $h$ to the sum $\sum
d_\lambda^p$ can be estimated as follows. By Theorem~\ref{thm:wav-hol}, the lower H\"older exponent $\lh(x)$ of a function at
$x$ is
\[
\lh(x)=\liminf_{j\to\infty \atop \lambda'\subset3\lambda_j(x)} \frac{\log d_{\lambda'}}{\log 2^{-j}},
\]
which allows us to write $d_\lambda\sim 2^{-h j}$. Moreover, the
number of these dyadic intervals should be about $2^{d(h)j}$, each
of volume $2^{-dj}$. Hence, the contribution is $2^{(d(h)-d-h
p)j}$. The dominating contribution is the one corresponding to the
value $h$ associated with the biggest exponent; by writing the
equality (\ref{eq:omega}) as $\sum d_\lambda^p\sim
2^{-\omega(p)j}$, one can expect the following relation,
$-\omega(p)=\sup_{h}\{d(h)-d-h p\}$. As $-\omega$ is a convex
function, if $d$ is concave, then $-\omega$ and $-d$ are convex
conjugate functions, so that $d(h)=\inf_p\{h p-\omega(p)+d\}$. Let
us remark that the preceding argument is far from being a
mathematical proof; see \cite{ajl:04} and \cite{jaf-nic:08} for a
comparison between the WLM and other multifractal formalisms.

Although it can be shown that formula (\ref{eq:smf:leg}) allows to
recover the spectrum of singularities under additional assumptions
(see \cite{jaf:96,jaf:97,aou-ben:02} for instance), the validity
does not hold  in complete generality. Indeed, the only result
valid in the general case is the following inequality
\cite{jaf:97,jaf:04},
\begin{equation}\label{eq:smf:ineg}
\underline{d}(h)\le \inf_{p\in\R_*} \{hp-\omega(p)+d\}.
\end{equation}

Once the lower spectrum of singularities has been introduced, the
upper spectrum of singularities can be defined in a totally
analogous way. A relation similar to the inequality
(\ref{eq:smf:ineg}) holds.

The weak multifractal formalism is defined as follows. \BDef Let
$f$ be a locally bounded function; its upper isoH\"older sets are
the sets
\[
\overline{E}_H=\{x:\uh(x)=H\}.
\]
The upper spectrum of singularities of $f$ is the function
\[
\overline{d}:
\R^+\cup \{\infty\}\to \R^+\cup \{-\infty\} \quad
H\mapsto \dim_\hau(\overline{E}_H).
\]
\EDef

The following theorem which can be found in \cite{ajvw:08} gives
an upper bound for the upper spectrum of singularities.
\BThe Let
$f$ a uniform H\"{o}lder function. The following inequality holds
\begin{equation}\label{eq:majospectr}
\overline{d}(h)\le \inf_{p\in\R^*} \{hp-\omega(p)+d\}.
\end{equation}
\EThe
\BDef Let $f$ a uniform H\"{o}lder function and $h>0$; if
(\ref{eq:majospectr}) is an equality, i.e.
\[
\forall h>0,\quad \overline{d}(h)= \inf_{p\in\R^*} \{hp-\omega(p)+d\}
\]
then function $f$ is said to obey the weak multifractal formalism.
\EDef
\section{Construction of a class of wavelet series obeying the weak multifractal
formalism}\label{sec4}
The aim of this Section is to exhibit a
class of multifractal functions for pointwise irregularity. This
question is in fact a non trivial one. A quite natural
approach to solve this problem is to consider multifractal
functions for the usual pointwise regularity.

Let us point out that if we want to define wavelet series
that are multifractal both for the strong and weak H\"{o}lderian regularity point of views, we
have to take into account that, except in the case where the lower and
upper exponents coincide, there is no wavelet criteria for
the pointwise irregularity.

In the same spirit as Barral and Seuret in \cite{bs:05}, we
will define wavelets series built from a multifractal measure
$\mu$ on $[0,1[^{d}$ in the following way,
\begin{equation}\label{EqDefWSeries}
F_{\mu}(x)=
\sum_{j\geq 0} \sum_{|\lambda|=2^{-j}, \atop \lambda\subset
[0,1]^{d}}2^{-j(s_{0}-\frac{d}{p_{0}})}\mu(\lambda)^{\frac{1}{p_{0}}}\psi_{\lambda}(x),
\end{equation}
$\forall x\in [0,1]^{d}$, where the wavelets $\psi^{i}$ belongs to the Schwartz class on
$\R^{d}$ with all moments vanishing. This class of examples also
proves that  upper and lower spectra may coincide: Under specific assumptions detailed in section~\ref{SecClassFct}, we can obtain a class of functions obeying both the strong and the weak multifractal formalisms.

We begin by recalling some basic facts about the multifractal analysis of
measures.
\subsection{Some results about multifractal analysis of
measures}\label{subseq:Some results about multifractal analysis of
measures}
Following Barral and Seuret, we adapt here the usual
multifractal formalism of \cite{bmp:92}. The main difference lies
in the definition of the isoH\"{o}lder sets, since we just need a
multifractal formalism
associated with a dyadic grid.\\
We first give some slightly modified versions of the usual
definitions of lower and upper exponents
of a given Borel measure $\mu$ at a point $x_{0}$.
For any $\sigma\in\{-1,0,1\}^d$ and any dyadic
cube
$\lambda=\prod_{\ell}[\frac{k_{\ell}}{2^{j}},\frac{k_{\ell}}{2^{j}+1}[$,
let us set
\[
\lambda^{\sigma}=\prod_{\ell}[\frac{k_{\ell}+\sigma_{\ell}}{2^{j}},
\frac{k_{\ell}+1+\sigma_{\ell}}{2^{j}}[ \quad
\mbox{and}\quad \mu^{\sigma}(\lambda)=\mu(\lambda^{\sigma}).
\]
We also define the quantities
\[
\underline{\alpha}_{\mu}^{\sigma}(x_{0})=\liminf_{j\to
\infty} \frac{\log \mu^\sigma(\lambda_{j}(x_0))}{-j\log(2)},\quad
\overline{\alpha}_{\mu}^{\sigma}(x_{0})=\limsup_{j\to
\infty} \frac{\log \mu^\sigma(\lambda_{j}(x_0))}{-j\log(2)},
\]
and, in case of existence,
\[
\alpha_{\mu}^{\sigma}(x_{0})=\lim_{j\to
\infty} \frac{\log \mu^\sigma(\lambda_{j}(x_0))}{-j\log(2)}.
\]
We will be concerned by the estimate of the Hausdorff dimension of
the following isoH\"{o}lder sets
\[
\widetilde{E}_{\alpha}(\mu)=\{x\in
[0,1[^{d},\,\alpha_{\mu}(x)=\min\limits_{\sigma\in\{-1,0,1\}^d}\{\alpha_{\mu}^{\sigma}(x)\}=\alpha\}.
\]
The mapping
\[
d_{\mu}:\alpha\geq 0\mapsto
\dim_{\hau}(\widetilde{E}_{\alpha}(\mu)),
\]
will be called the multifractal spectrum of the Borel measure
$\mu$. Recall that in the framework of \cite{bmp:92}, the
following isoH\"{o}lder sets are used
\[
E_{\alpha}(\mu)=\{x\in [0,1[^{d},\,\lim_{j\to
+\infty} \frac{\log \mu(\lambda_j(x))}{-j\log(2)}=\alpha\}.
\]
Unfortunately, these isoH\"{o}lder sets are not adapted to the
study of the pointwise regularity of wavelet series $F_{\mu}$.
Indeed, starting from $\lim_{j} \log
\mu(\lambda_j(x))/-j\log(2)$, we cannot deduce the value of
the upper pointwise H\"{o}lder exponent of the function $F_{\mu}$
at $x$ using wavelet criteria.

We now recall well known results about upper bound of the upper
multifractal
spectrum, which can be found in \cite{bmp:92}. For any $q$ in $\R$, set
\[
\tau(q)=\liminf_{j\to
+\infty} \frac{1}{\log |\lambda|}\log \sum\limits_{|\lambda|=2^{-j}}^{*}\mu(\lambda)^{q},
\]
where $\sum\limits^{*}$ means that the sum is taken over those $\lambda$
such that $\mu(\lambda)>0$.\\
As usual, $\tau^{*}$ denotes the Legendre transform of the function
$\tau$, that is
\[
\forall \alpha\geq 0,\quad
\tau^{*}(\alpha)=\inf_{q\in\R} \{\alpha q-\tau(q)\}.
\]
Remark that, since $\alpha>0$,
\[
\widetilde{E}_{\alpha}(\mu)\subset E_{\alpha}(\mu).
\]
Using this inclusion, an upper bound for the multifractal spectrum
of any Borel measure can be obtained from \cite{bmp:92}: \BPro Let
$\alpha\geq 0$ and $\mu$ a Borel measure. One has
\[
\dim_{\hau}(\widetilde{E}_{\alpha}(\mu))\leq\tau^{*}(\alpha).
\]
Moreover, if $\tau^{*}(\alpha)<0$ then
$\widetilde{E}_{\alpha}(\mu)=\emptyset$. \EPro \BDef Let
$\alpha_0\geq 0$. One says that the Borel measure $\mu$ obeys the
multifractal formalism at $\alpha=\alpha_0$ for the sets
$\widetilde{E}_{\alpha}(\mu)$
 if
$\dim_{\hau}(\widetilde{E}_{\alpha_0}(\mu))=\tau^{*}(\alpha_0)$.
\EDef

\subsection{Wavelet series and multifractal measures}
We want to define a wavelet series of the form
(\ref{EqDefWSeries}) obeying the weak multifractal formalism for
functions. First, we give an explicit relationship between the wavelet series
$F_{\mu}$ and the measure $\mu$ from the multifractal point of view.
\subsubsection{A transference theorem}
\BThe\label{TheoBarSeur}
Let $\mu$ a Borel measure
  and $s_{0},p_{0}$
 two positive real numbers.
Let $F_{\mu}$ be the wavelet series defined by equality
(\ref{EqDefWSeries}). If the measure $\mu$ obeys the multifractal
formalism at $\alpha_0\geq 0$ for the sets
$\widetilde{E}_{\alpha}(\mu)$, then $F_{\mu}$ obeys both the
strong and weak multifractal formalisms at
\[
H=s_{0}-\frac{d}{p_{0}}+ \frac{\alpha_0 d}{p_{0}},
\]
and
\[
\underline{d}(H)=\overline{d}(H)=d_{\mu}(\alpha_0).
\]
\EThe
\BRem Let us notice, as in \cite{bs:05}, that if $x_{0}\not\in
supp(\mu)$, there exists some $j_{0}$ such that,
\[
\forall j\geq j_{0},\quad d_{j}(x_{0})=0.
\]
Thus, in this special case,
$\overline{h}_{F_{\mu}}(x_{0})=+\infty$.
\ERem
\BProof
The proof mimics the one of Theorem $1$ of \cite{bs:05}. It relies
on the following Lemma:
\BLem For any $\alpha\geq 0$, the
following inclusion holds:
\[
\widetilde{E}_{\alpha}(\mu)\subset
\underline{E}_{H}(F_{\mu})\cap\overline{E}_{H}(F_{\mu}),
\]
where
$H=s_{0}-\displaystyle\frac{d}{p_{0}}+\displaystyle\frac{\alpha
d}{p_{0}}$.\ELem
\BProof
Remark that, for all $\lambda$
\[
d_{\lambda}=2^{-j(s_{0}-\frac{d}{p_{0}})}\mu(\lambda)^{\frac{1}{p_{0}}}.
\]
where $j=-\log(|\lambda|)/\log(2)$.\\
For any given $x_{0}$,
\[
d_{j}(x_{0})=2^{-j(s_{0}-\frac{d}{p_{0}})}\max_{\sigma\in\{-1,0,1\}^d} \{\mu^{\sigma}(\lambda_{j}(x_0))\}^{\frac{1}{p_{0}}}.
\]
Assume that $x_0\in \widetilde{E}_{\alpha}(\mu)$. Since for any
$\sigma\in\{-1,0,1\}^d$, $\alpha^{\sigma}(x_{0})\geq \alpha$, for
all $\varepsilon>0$ there exists an integer
$j_{0}(\varepsilon,\sigma)$ such that
\[
\forall j\geq j_{0}(\varepsilon,\sigma),\quad
\mu^{\sigma}(\lambda_{j}(x_{0}))\leq
2^{-j\frac{d}{p_0}(\alpha-\varepsilon)}.
\]
Hence,
\[
\forall j\geq \max_{\sigma} \{j_{0}(\varepsilon,\sigma)\},\quad
d_{j}(x_{0})\leq
2^{-j(s_{0}-\frac{d}{p_{0}})}2^{-j\frac{d}{p_0}(\alpha-\varepsilon)}\leq
2^{-j(H-\varepsilon)}
\]
and
\[
\liminf_{j\to
+\infty} \frac{\log d_{j}(x_{0})}{\log 2^{-j}}\geq
H.
\]
Furthermore, since for some $\sigma_0\in\{-1,0,1\}^d$
we have $\alpha^{\sigma_0}(x_{0})\leq \alpha$, for all $\varepsilon>0$
there exists an integer $j_{0}(\varepsilon)$ such that
\[
\forall j\geq j_{0}(\varepsilon),\quad
\mu^{\sigma_0}(\lambda_{j}(x_{0}))\geq
2^{-jd(\alpha+\varepsilon)}.
\]
Then,
\[
\forall j\geq j_{0}(\varepsilon),\quad
 d_{j}(x_{0})\geq
2^{-j(H+\varepsilon)}
\]
and
\[
\limsup_{j\to
+\infty} \frac{\log d_{j}(x_{0})}{\log 2^{-j}}\leq
H.
\]
Hence, using Corollary \ref{cor:specialcase},
\[
\widetilde{E}_{\alpha}(\mu)\subset \{x\in
[0,1]^d,\,\underline{h}_{F_{\mu}}(x)=\overline{h}_{F_{\mu}}(x)=H\}=
\underline{E}_{H}(F_{\mu})\cap\overline{E}_{H}(F_{\mu}).
\]
\EProof

Since $\mu$ obeys the
multifractal formalism at $\alpha=\alpha_0$ for the sets
$\widetilde{E}_{\alpha}(\mu)$,
\[
\dim_{\hau}(\widetilde{E}_{\alpha_0}(\mu))=\tau^{*}(\alpha_0)\leq
\dim_{\hau}(\overline{E}_{H}(F_{\mu}))
\]
with $H=s_{0}-\displaystyle\frac{d}{p_{0}}+\displaystyle\frac{\alpha_0 d}{p_{0}}$.\\
Remark then that for any $p\in\R$,
\[
\omega_f^{*}(p)=
p(s_{0}-\frac{d}{p_{0}})-\tau(\frac{p}{p_{0}}).
\]
Since for any $0<H<\infty$ and any locally bounded function one
has
\[
\dim_{\hau}(\overline{E}_{H}(F_{\mu}))\leq
\inf_{p\in\R}\{pH-\omega_{f}^{*}(p)+d\}=\tau^{*}(\alpha_0)=\dim_{\hau}(\widetilde{E}_{\alpha_0}(\mu)),
\]
one can conclude that $F_\mu$ obeys the weak multifractal formalism. A
similar approach proves that $F_\mu$ also obeys the strong
multifractal formalism.
\EProof

\subsection{A class of wavelet series obeying both the strong and the weak
multifractal formalisms}\label{SecClassFct} The aim of this
section is to exhibit a class of multifractal measures obeying the
multifractal formalism at any $\alpha\geq 0$ for the sets
$\widetilde{E}_{\alpha}(\mu)$, yielding an example of wavelet
series satisfying both the strong and weak multifractal
formalisms. To this end we first give some examples of multifractal measures obeying the multifractal formalism for sets $\widetilde{E}_{\alpha}$ using Theorem $2$ of \cite{bs:05} for sets $\widetilde{E}_{\alpha}$. Indeed, even if we consider slightly different iso--H\"{o}lder sets Theorem 2 still holds : the proof is exactly the same that this of \cite{bs:05}.\\
We give two canonical examples
of measures satisfying the conditions above.
\subsubsection{Quasi-Bernoulli measures} Let $b$ an integer larger than
$2$. Let us recall that a Borel positive measure on $[0,1]^d$
$\mu$ is said quasi-Bernoulli if for some $C>0$ and for any
$v,w\in (\mathcal{A}^d)^{*}$,
\[
\displaystyle\frac{1}{C}\mu(I_{v})\mu(I_{w})\leq \mu(I_{vw})\leq
C\mu(I_{v})\mu(I_{w}).
\]
A classical example of quasi-Bernoulli measures is the well-known
example of multinomial measures: \BEx Let $b$ an integer larger
than $2$ and let $(m_{0},\cdots,m_{b-1})\in (0,1)^{b}$ such that,
\[
\sum\limits_{i=0}^{b-1}m_{i}=1;
\]
we can construct a sequence of probability measures
$(\mu_{n})_{n\in\N}$ on $[0,1)^{d}$ as follows. For any integer
$n$, define a probability measure $\mu_{n}$ on $[0,1)^d$ such that
for any $w\in\mathcal{A}^{nd}$,
\[
\mu_{n}(I_{w})=\prod\limits_{\ell=1}^{nd}m_{w_{\ell}}.
\]
This sequence has a weak limit $\mu$ called multinomial measure of
base $b$ with weight $(m_{0},\cdots,m_{b-1})$.\\
By construction, any multinomial measure is quasi-Bernoulli.
\EEx
In the following, we consider only continuous quasi-Bernoulli
measure, that is without atom. Recall that any
continuous quasi-Bernoulli measure satisfies the assumptions of
Theorem 2 of \cite{bs:05} and thus obeys the multifractal formalism at
any $\alpha>0$ for the sets $\widetilde{E}_\alpha$.\\

Hence we have the following
Theorem.
\BThe
If $\mu$ is a continuous quasi-Bernoulli measure,
then the wavelet series $F_\mu$ defined by (\ref{EqDefWSeries})
obeys both the strong and the weak multifractal formalisms at any
$H>s_0- d/p_0$.
\EThe

\subsubsection{The case of $b$-adic random multiplicatives cascades}
Let $b$ an integer larger than $2$ and $d=1$. Canonical random
cascades were introduced by Mandelbrot in \cite{man:74} and their
multifractal properties have been widely studied, mainly in the
setting of $b$-adic grid (see e.g.\
\cite{kah:91,hw:92,ck:92,molc:96,bar:99,bar:00}). We first recall
the construction of these measures. Let $W$ a non negative random
variable, not almost surely constant, satisfying $\E(W)=1/b$. We
thus consider $(W_w)_{w\in \mathcal{A}^{*}}$ a sequence of
independent copies of $W$ and $\mu_n$ the random measure whose
density with respect to the Lebesgue measure on any dyadic
interval is constant and equals
\[
b^n W_{w_1}\cdots W_{w_1\cdots w_n}.
\]
Almost surely, this sequence of measures converges weakly to a
measure $\mu$ as $n$ goes to infinity. Recall that if $\mu$ is a
$b$-adic random multiplicative cascade, then almost surely on $J$,
$\mu$ satisfies the assumptions of Theorem 2 of \cite{bs:05} and thus
obeys the multifractal formalism for any $q\in J$ at
$\alpha=\widetilde{\tau}'(q)$ for the sets $\widetilde{E}_\alpha$.\\

Then similarly, to the case of quasi-Bernoulli measure, we have the following result,
\BThe
Let $W$ be an almost surely positive random variable. Let $\mu$
a $b$-adic random multiplicative cascade such that
$\widetilde{\tau}'(1)=-1-\log_b(\E(W))>0$. Then the
wavelet series $F_\mu$ defined by (\ref{EqDefWSeries}) obeys
almost surely both the strong and the weak multifractal formalisms at
any $H>0$.
\EThe

\section{A multifractal function whose lower and upper spectra of singularities differ}
Although similar results hold for both the strong and the weak
multifractal formalisms, there is no direct relation between
$\underline{d}$ and $\overline{d}$. We introduce here a
function defined as $p$-adic Davenport series whose upper multifractal
spectrum is reduced to two single points, while
its lower multifractal spectrum is linear on an interval.

A $p$-adic Davenport series ($p\ge 2$) is a series of the form
\[
f(x)=\sum_{j=0}^\infty a_j \fp{p^j t},
\]
where $\fp{x}$ is the sawtooth function
\[
\fp{x}=x-[x]-\frac{1}{2}.
\]
We will assume here that $(a_j)_j\in l^1$, so that the series is
normally convergent. The function $f$ is thus continuous at every
non $p$-adic rational number and has left and right limit at every
$p$-adic rational $kp^{-l}$ ($k\wedge p=1$) with a jump of
amplitude $\sum_{m\ge l} a_m$. Recent results on Davenport series can be found in \cite{jaf:04:bis}. Let $\beta>1$; the functions $f_\beta$ we will study is defined by
\[
f_\beta(x)=\sum_{l\in\N} \frac{\fp{2^lx}}{2^{l\beta}}.
\]

\subsection{The lower spectrum of singularities of $f_\beta$}
The functions $f_\beta$ are derived from the famous
L\'evy's function (which can be seen as a special case, where
$\beta=1$). The properties of the lower spectrum of singularities
of this function have already been investigated in \cite{jaf:97};
Propositions~\ref{pro:lh fbeta} and \ref{pro:ss fbeta} together
can be seen as a generalization of Proposition~4 of \cite{jaf:97}.
Proposition~12 of \cite{jaf:04:bis} implies that $\underline{d}$ is
linear on $[0,\beta]$.

To determine explicitly the lower isoH\"older sets of $f_\beta$,
we will use the following notations. Let $p\in\N$, $p>1$; for a
sequence of integers $(x_l)_{l\in\N}$ satisfying $0\le x_j<p$, we
will write
\begin{equation}\label{eq:expansion}
(0;x_1,\ldots,x_l,\ldots)_p
\end{equation}
to denote one expansion in basis $p$ of the real number
\[
x=\sum_{l\in\N} \frac{x_l}{p^l}.
\]
If there is no $k$ such that $x_l=p-1$ for all $l\ge k$,
(\ref{eq:expansion}) is the proper expansion of $x$ in basis $p$.
If $(0;x_1,\ldots)_p$ is the proper expansion of $x$, we define
\[
\theta_p(x)=\inf\{l:x_l\not=0\}-1.
\]
Let $\delta(k)=\sup\{l:\forall l'\le l, x_{k+l'}=x_{k}\}$ and let $(m_l)_{l\in\N}$ be the sequence defined recursively,
$m_1=\inf\{l:x_l=0 \mbox{ or } x_l=p-1\}$,
$m_k=\inf\{l\ge m_{k-1}+\delta(m_{k-1}): x_l=0\mbox{ or } x_l=p-1\}$
($k>1$). One also defines the sequence $(\delta_k)_{k\in\N}$ by
$\delta_k=\delta(m_k)$. Finally, $\rho_p(x)=\limsup_{k\to\infty}
\delta_k/m_k$; if $x$ is a $p$-adic rational, one sets
$\rho_p(x)=\infty$. The number $\rho_p(x)$ defines, in some way,
the rate of approximation of the number $x$ by $p$-adic rationals, since we have the following obvious result.
\BPro
If $x$ is not a $p$-adic rational, the equation (depending on $k$ and $l$)
\[
|x-\frac{k}{p^l}|\le (\frac{1}{p^l})^\phi \qquad (k\wedge p=1)
\]
has an infinity of solutions if  and only if $\phi\le \rho_p(x)+1$.
\EPro
We will denote by $\phi(x)$ the critical exponent $\phi(x)=\rho_2(x)+1$.
The lower H\"older exponents of $f_\beta$ only depend on $\phi$.
\BPro\label{pro:lh fbeta} The lower H\"older exponents of
$f_\beta$ are given by
\[
\lh(x)=\frac{\beta}{\phi(x)}.
\]
\EPro
\BProof As a corollary of Theorem~21 of
\cite{jaf-nic:08}, we have the following equalities: if $x$ is not
a dyadic rational,
\begin{equation}\label{eq:levy-h}
\lh(x)=\liminf_{j\to \infty} \frac{-\beta j}{\log_2
\dist(x,2^{-j}\Z)};
\end{equation}
otherwise, $\lh(x)=0$. We can suppose that $x\in(0,1)$ is not a
dyadic rational. For a given $j\in\N$, let
$\varepsilon_j=\dist(x,2^{-j}\Z)$. One has
$\theta_2(\varepsilon_j)=j+1+\delta(j+1)$ and thus $\varepsilon_j\sim
2^{-(j+\delta(j+1)+1)}$. Then (\ref{eq:levy-h}) can be rewritten
\[
\lh(x)=\liminf_{j\to\infty} \frac{\beta
j}{j+1+\delta(j+1)}=\frac{\beta}{1+\rho_2(x)}.
\]
\EProof The lower isoH\"older sets are now characterized. \BCor The
lower isoH\"older sets of the function $f_\beta$ are the sets
\[
\underline{E}_H=\{x: \phi(x)=\frac{\beta}{H}\}\quad
(0<H\le\beta).
\]
The set $\underline{E}_0$ is the set of the dyadic rationals.
\ECor To conclude this study on the strong H\"older regularity, we
have the following result.
\BPro\label{pro:ss fbeta}
The lower spectrum of singularities of $f_\beta$ is
\[
\underline{d}(h)=\left\{\begin{tabular}{ll}
$\hphantom{-}\frac{h}{\beta}$ & if $h\in [0,\beta]$ \\
$-\infty$ & otherwise
\end{tabular}\right..
\]
\EPro
\BProof The main idea is the same as in Proposition~4 of
\cite{jaf:97}. If $\alpha\ge 1/\beta$, let
\[
F_\alpha=\limsup_{j\to\infty} \bigcup_k
[k2^{-j}-2^{-j\alpha\beta},k2^{-j}+2^{-j\alpha\beta}].
\]
Using (\ref{eq:levy-h}), $\lh(x)=H$ means
\begin{equation}\label{eq:dem-equiv}
x\in
\bigcap_{\gamma>H}F_{1/\gamma}-\bigcup_{\gamma<H}F_{1/\gamma}.
\end{equation}
Clearly, $\dim_\hau(F_\alpha)\le 1/\alpha\beta$; let us show that
the converse inequality holds.

Let $(j_l)_{l\in\N}$ be a sequence satisfying $j_l=2^{j_{l-1}}$,
let
\[
I_k(l)=[k2^{-j_l}-2^{-j_l
\alpha\beta},k2^{-j_l}+2^{-j_l\alpha\beta}]
\]
and
\[
G_\alpha=\bigcap_l \bigcup_k I_k(l).
\]
A probability measure $\mu$ supported by $G_\alpha$ can be
obtained as follows. If $l=1$, we put on each interval $I_k(1)$
the same mass $2^{-j_1}$. If each of these intervals contains $n$
intervals of type $I_k(2)$, on each of these intervals, we put the
measure $2^{-j_1}/n$. This construction can be iterated to obtain,
at the limit, a probability measure $\mu$ supported by $G_\alpha$.
One easily checks that
\[
\mu([x-h,x+h])\le C h^{1/\alpha\beta}\quad \forall x\in G_\alpha.
\]
Moreover, Proposition~4.9 of \cite{fal:90} implies that
\[
\hau^{1/\alpha\beta}(G_\alpha)>0
\]
and thus, since $G_\alpha\subset F_\alpha$,
\[
\dim_\hau(F_\alpha)=1/\alpha\beta,
\]
which, thanks to (\ref{eq:dem-equiv}), is sufficient to conclude.
\EProof

\subsection{The upper spectrum of singularities of $f_\beta$}
We show here that from the weak H\"older
regularity point of view, the function $f_\beta$ only displays two
kinds of singularities: it is discontinuous at dyadic rationals
and has an upper H\"older exponent equal to $\beta$ at non dyadic
rationals.

 Let
\[
\Omega_\alpha=\liminf_{j\to\infty} \{x\in\R: \exists k\in\Z \mbox{
such that } |x-\frac{k}{2^j}|\le 2^{-\alpha j}\}.
\]
We have the following relation between the sets $\Omega_\alpha$
and $\uh(x_0)$. \BPro If $\alpha>1$, then
\[
x_0\notin\Omega_\alpha \Rightarrow \uh(x_0)\ge
\frac{\beta}{\alpha}.
\]
\EPro
\BProof Let $\varepsilon>0$ be a given real number. We want to
prove that for any $C>0$, there exists a strictly decreasing
sequence $(r_n)_n$ of real positive numbers, such that
\[
\sup_{h>0} \|\Delta_h^{\! [\beta/\alpha]+1} f_\beta
(x)\|_{L^\infty(B_h(x_0,r_n))}\le C r_n^{\beta/\alpha}.
\]
One has
\begin{equation}\label{eq:dem-delta}
\Delta_h^{\! [\beta/\alpha]+1} f_\beta
(x)=\sum_{l\in\N_0}\frac{1}{2^{l\beta}} \sum_{m=0}^M (-1)^m
\bin{M}{m} \fp{2^l(x+mh)}.
\end{equation}
If $x_0\notin \Omega_\alpha$, then there exists a strictly
increasing sequence of integers $(j_n)_n$ such that
\[
|x_0-\frac{k}{2^{j_n}}|\ge 2^{-\alpha j_n},\quad \forall n\in\N,\
k\in\Z.
\]
Let $r_n=2^{-\alpha j_n}$; the interval
$[x,x+([\beta/\alpha]+1)h]\subset B(x_0,r_n)$ does not contain any
dyadic rational of the form $k2^{-l}$, with $l\le j_n-1$. This
implies
\[
\sum_{l=0}^{j_n-1} \frac{1}{2^{l\beta}} \sum_{m=0}^M (-1)^m
\bin{M}{m} \fp{2^l(x+mh)}=0.
\]
Relation (\ref{eq:dem-delta}) leads to the following inequality,
\begin{eqnarray*}
|\Delta_h^{\! [\beta/\alpha]+1} f_\beta
(x)|&=&|\sum_{l=j_n}^\infty\frac{1}{2^{l\beta}}
\sum_{m=0}^M (-1)^m \bin{M}{m} \fp{2^l(x+mh)}| \\
&\le& C' \sum_{l=j_n}^\infty \frac{1}{2^{l\beta}} \le C'
r_n^{\beta/\alpha}.
\end{eqnarray*}
Let now $n_0$ be an integer such that $C' r_n^{\beta/\alpha}\le C
r_n^{\beta/\alpha-\varepsilon}$ for all $n\ge n_0$. We have
$\uh(x_0)\ge\beta/\alpha-\varepsilon$, which is sufficient to
conclude.
\EProof
The sets $\Omega_\alpha$ are explicitly known
whenever $\alpha\ge 1$. \BPro If $\alpha>1$, then
\[
\Omega_\alpha=\{\frac{k}{2^j}:(k,j)\in \Z\times \N\}.
\]
Moreover, $\Omega_1=\R$.
\EPro
\BProof The case $\alpha=1$ is
trivial. Let $\alpha>1$ and let $x\in (0,1)$ be a non dyadic
rational. If $x=(0;x_1,\ldots)_2$, one has, for $j$ sufficiently
large,
\[
\min_{k\in\Z} |x-\frac{k}{2^j}| \ge 2^{-n(j)+1},
\]
where $n(j)$ is the first index greater or equal to $j$ such that
$x_{n(j)}=1$. We then have, for $j=n(j)-1$ sufficiently large,
\[
\min_{k\in\Z} |x-\frac{k}{2^j}| \ge 2^{-j}> 2^{-\alpha j}.
\]
Therefore, $x\notin \Omega_\alpha$.
\EProof
Thus, we have a lower
bound for the upper H\"older exponent. \BCor\label{cor:lower} If
$x_0$ is not a dyadic rational, then
\[
\uh(x_0)\ge \beta.
\]
\ECor Let us now prove the converse inequality. The following
proposition is similar to Lemma 1 of \cite{jaf:97}.
\BPro\label{pro:upper} Let $f$ be a function defined on $\R$,
continuous everywhere except on a dense countable set of points
and admitting a left and a right limit at every point. Let also
$x_0\in \R$ be a point of continuity of $f$ and $(r_n)_n$ a
sequence of points of discontinuity converging to $x_0$. Finally,
let $s_n$ ($n\in\N$) be the jump of $f$ at $r_n$. If there exists
a strictly increasing function $\psi$ satisfying $\psi(0)=0$ such
that
\[
\exists r_0: \Big( \forall r\le r_0, \exists r_n:
\big(|r_n-x_0|\le r,\ |s_n|\ge\psi(r)\big)\Big),
\]
then
\[
\uh(x_0)\le \limsup_{r\to 0} \frac{\log \psi(r)}{\log r}.
\]
\EPro
\BProof Let $\alpha>1$; if $f\in C^\alpha_w(x_0)$, for any
$C>0$ there exists a sequence $(t_n)_n$ such that
\[
\sup_{h>0} \|\Delta_h^{\!
[\alpha]+1}f(x)\|_{L^\infty(B_h(x_0,t_n))} \le C t_n^\alpha.
\]
Let $(r_n)_n$ a sequence such that, for any integer $n$
sufficiently large,
\[
|r_n-x_0|<t_,\quad |s_n|\ge \psi(t_n).
\]
Let $F_n$ be the function defined on $\R$ by
\[
F_n(h)=(f(r_n+h),\ldots,f(r_n+([\alpha]+1)h)).
\]
Since $F_n$ has only a countable set of discontinuities, one can
find $h$ arbitrarily close to zero such that
$[r_n,r_n+([\alpha]+1)h]\subset B(x_0,t_n)$ and such that $F_n$ is
continuous at $h$. Therefore,
\begin{eqnarray*}
\psi(t_n)&\le& s_n=|f(r_n^+)- f(r_n^-)| \\
&\le& |\sum_{k=1}^{[\alpha]+1} (-1)^k \bin{[\alpha]+1}{k} f(r_n+kh)+ f(r_n^-)|\\
&& +
|\sum_{k=1}^{[\alpha]+1} (-1)^k \bin{[\alpha]+1}{k} f(r_n+kh)+f(r_n^+)| \\
&\le& 2 (|\Delta_h^{\! [\alpha]+1} f(r_n^-)| + |\Delta_h^{\!
[\alpha]+1} f(r_n^+)|)\le 4C t_n^\alpha.
\end{eqnarray*}
This inequality implies
\[
\frac{\log\psi(t_n)}{\log t_n} \ge \alpha+ \frac{\log 4C}{\log
t_n}.
\]
As $t_n$ tends to zero,
\[
\limsup_{r\to 0} \frac{\log\psi(r)}{\log r}\ge\alpha,
\]
and the result follows.
\EProof
Since $f_\beta$ is continuous
except at dyadic rationals and since $\Omega_1=\R$, Proposition~\ref{pro:upper} and Corollary~\ref{cor:lower} imply the following
result.
\BThe
If $x_0$ is not a dyadic rational,
\[
\uh(x_0)=\beta,
\]
if $x_0$ is a dyadic rational, $\uh(x_0)=0$.
\EThe

\paragraph{Acknowledgment.} The authors would like to thank A. Ayache, Y. Heurteaux and S. Jaffard for several stimulating discussions.

\end{document}